\documentclass[12pt]{article}
\usepackage{mitpress}

\newtheorem{theorem}{Theorem}

\newtheorem{corollary}[theorem]{Corollary}

\newtheorem{definition}[theorem]{Definition}
\newtheorem{example}[theorem]{Example}

\newtheorem{lemma}[theorem]{Lemma}

\newtheorem{remark}[theorem]{Remark}

\newenvironment{proof}[1][Proof]{\textbf{#1.} }{\ \rule{0.5em}{0.5em}}
\newdimen\dummy
\dummy=\oddsidemargin
\input{tcilatex}

\begin{document}

\title{Homeomorphisms of $\overline{U}\times R$ and rotation number}
\author{Paul Fabel \\
Drawer MA\\
Mississippi State University\\
fabel@ra.msstate.edu}
\maketitle

\begin{abstract}
Suppose $U\subset R^{2}$ is bounded open and contractible and $H:\overline{U}%
\times R\rightarrow \overline{U}\times R$ is a homeomorphism leaving
invariant $U\times R$. If $\partial U$ is locally connected and not a simple
closed curve $H$ induces a homeomorphism of the solid cylinder leaving
invariant sufficiently many vertical lines to determine a rotation number.
If $\partial U$ is not locally connected $H$ admits a natural notion of
rotation number despite a general absence of an induced homeomorphism of the
solid cylinder.
\end{abstract}

\section{Introduction}

The following theorem, credited to Ursell and Young \cite{ursell} is central
to applications (\cite{alli},\cite{barge},\cite{barge1},\cite{brech},\cite
{cart},\cite{mayer}) of prime end theory to dynamical systems:

\begin{theorem}
\label{prmend}Suppose $U\subset R^{2}$ is bounded open and contractible, $%
\psi :U\rightarrow int(D^{2})$ is conformal, and $h:\overline{U}\rightarrow 
\overline{U}$ is a homeomorphism such that $h(U)=U.$ Then $\psi h\psi ^{-1}$
can be extended to a homeomorphism of $D^{2}.$
\end{theorem}

In particular if $h$ is orientation preserving then $\overline{\psi h\psi
^{-1}}_{S^{1}}:\partial D^{2}\rightarrow \partial D^{2}$ determines a
rotation number, measuring in a sense the average rotation by $h$ of $%
\partial U$ about $U$.

Specific examples of higher dimensional versions of Theorem \ref{prmend} are
in short supply despite a general criteria established in \cite{brech}.

For example if $\partial U$ is not locally connected, there is no
corresponding version of Theorem \ref{prmend} for domains $U\times R\subset
R^{3}$ with $\Psi :U\times R\rightarrow int(D^{2})\times R$ defined via $%
\psi (u,t)=(\psi (u),t)).$ Example \ref{newpr} exhibits a homeomorphism $H:%
\overline{U}\times R\rightarrow \overline{U}\times R$ such that $H(U\times
R)=U\times R$ but $\Psi H\Psi ^{-1}$ \textbf{cannot} be extended to a
homeomorphism of $D^{2}\times R$.

The main result of this paper (Theorem \ref{mainthm}) salvages a notion of
rotation number for such homeomorphisms $H:\overline{U}\times R\rightarrow 
\overline{U}\times R.$ The basic idea is that $H$ preserves the circular
order of a certain collection of sets, each of which can be understood as
the product of $R$ with an interval of accessible prime ends of $U.$ If $%
\partial U$ is not locally connected this provides enough information to
determine a homeomorphism $g:S^{1}\rightarrow S^{1}$ whose rotation number
we declare to be that of $H.$

On the other hand a useful `3 dimensional prime end theory'\ exists if $%
\partial U$ is locally connected. Theorem \ref{pelc} shows $\Psi $ induces a
homeomorphism of the (3 cell) two point compactification of $D^{2}\times R.$
Moreover if $\partial U$ is not a simple closed curve of the cutpoints of $%
\partial U$ help to determine a discrete collection of invariant boundary
lines which in turn determine a rotation number.

Both notions of rotation number are invariant under topological conjugacy,
and agree with the usual rotation number of $h$ in the special case $%
H(u,t)=(h(u),t)$ where $h:\overline{U}\rightarrow \overline{U}$ is a
homeomorphism such that $h(U)=U.$

\section{Preliminaries}

Suppose throughout this paper $U\subset R^{2}$ is bounded open and
contractible, $\partial U$ denotes $\overline{U}\backslash U,D^{2}\subset
R^{2}$ denotes the closed unit disk and $\psi :U\rightarrow int(D^{2})$ is
conformal. All function spaces will have the compact open topology.

Define $\Psi :U\times R\rightarrow int(D^{2})\times R$ via $\Psi (u,t)=(\psi
(u),t).$

Let $D^{3}=\{(x,y,z)\in R^{3}|\sqrt{x^{2}+y^{2}+z^{2}}\leq 1\},$ the
standard 3 cell.

If $J\subset S^{1}$ is connected let $int(J)=J$ if $J=\{x\}.$ Otherwise let $%
int(J)$ denote the union of all open intervals contained in $J.$

If $Y$ is a topological space attach two points $\{\infty \}$ and $\{-\infty
\}$ to $Y\times R$ creating $(Y\times J)\cup \{\infty ,-\infty \}$
topologized such that $(y_{n},t_{n})\rightarrow \{-\infty \}$ iff $%
t_{n}\rightarrow \infty $ and $(y_{n},t_{n})\rightarrow \{\infty \}$ iff $%
t_{n}\rightarrow -\infty .$

In similar fashion we attach two points $\{\infty \}$ and $\{-\infty \}$ to $%
Y\times (-1,1)$ and create a new space $\overline{Y\times (-1,1)}$
topologized such that $(y_{n},t_{n})\rightarrow \{-\infty \}$ iff $%
t_{n}\rightarrow -1$ and $(y_{n},t_{n})\rightarrow \{\infty \}$ iff $%
t_{n}\rightarrow 1.$

\section{The map $\protect\phi :X^{\ast }\rightarrow D^{2}$}

The following procedure creates a complete metric space $X^{\ast }$ whose
underlying set can be seen as the union of $U$ and the accessible prime ends
of $U$.

Define a metric $d^{\ast }:U\times U\rightarrow R$ such that $d^{\ast
}(x,y)<\varepsilon $ iff there exists a map $f:[0,1]\rightarrow U$ such that 
$f(0)=x,f(1)=y$ and $\forall s,t$ $\left| f(t)-f(s)\right| <\varepsilon .$
Let $(X^{\ast },d^{\ast })$ denote the metric completion of $(U,d^{\ast }).$
Let $\partial X^{\ast }=X^{\ast }\backslash U^{\ast }.$

\begin{lemma}
\label{smallmain}There exists a map $\phi :X^{\ast }\rightarrow D^{2}$ such
that $\phi $ is uniformly continuous, one to one and $\phi (U^{\ast
})=int(D^{2}).$
\end{lemma}

\begin{proof}
Let $\overline{id}:X^{\ast }\rightarrow \overline{U}$ denote the unique
extension of the uniformly continuous identity map $id:U^{\ast }\rightarrow
U.$ Uniform continuity of $\psi (id)$ is essentially a consequence of
Theorem \ref{prmend} and is proved in Theorem 3.1 of \cite{fabel}.Define $%
\phi :X^{\ast }\rightarrow D^{2}$ to be the unique continuous extension of $%
\psi (id).$ Suppose $x\neq y$ and $\{x,y\}\subset X^{\ast }.$ If $\overline{%
id}(x)\neq \overline{id}(y)$ it follows from Proposition 2.14 \cite{pom}
that $\psi (\overline{id}(x))\neq \psi \overline{id}(y).$ If $\overline{id}%
(x)=\overline{id}(y)$ then $\{\overline{id}(x),\overline{id}(y)\}\subset
\partial U$ and $\{x,y\}\subset \partial X^{\ast }.$ Construct a closed
topological disk $E\subset \overline{U}$ such that $\overline{id}(x)\in
\partial E$ and $E\backslash \{\overline{id}(x)\}\subset U.$ If $d^{\ast
}(x,y)\neq 0$ then $int(E)\cap \partial U\neq \emptyset $ and $\overline{id}%
(x)$ and $\overline{id}(y)$ determine distinct prime ends. Let $z\in
\partial E\backslash \{\overline{id}(x)\}.$ By Theorem 2.15 \cite{pom} $\psi 
$ determines a bijection between the prime ends of $U$ and $\partial D^{2}.$
In particular $\partial E\backslash \{z\}$ determines distinct endcuts which
map under $\psi $ to distinct points of $\partial D^{2}.$ Hence $\phi $ is
one to one.
\end{proof}

\begin{remark}
The injective map $\phi :X^{\ast }\rightarrow D^{2}$ need not be an
embedding. For example if $U\subset R^{2}$ is the region bounded by a
`Warsaw circle'. The canonical map from $\partial X^{\ast }\rightarrow
\partial D^{2}$ is a continuous bijection but not a homeomorphism.
\end{remark}

\section{Failure of homeomorphism extension}

\begin{example}
\label{newpr}Suppose $U\subset R^{2}$ is the interior of the standard Warsaw
circle. There there exists $H\in G(U\times R)$ such that $\Psi H\Psi
^{-1}:int(D^{2})\times R\rightarrow int(D^{2})\times R$ is not extendable to
a homeomorphism of $D^{2}\times R.$
\end{example}

\begin{proof}
Let $U\subset R^{2}$ be the interior of a standardly embedded Warsaw circle
such that the impression of the bad prime end ( the closed interval `limit
bar') is precisely $\partial U\cap (\{0\}\times R).$ Define a homeomorphism $%
H:U\times R\rightarrow U\times R$ via $H(x,y,t)=(x,y,y+t).$ Note $\Psi H\Psi
^{-1}$ is not continuously extendable near the bad prime end since $H$ is
not level preserving on the impression of the bad prime end.
\end{proof}

\section{Homeomorphisms of $S^{1}$ and rotation number}

The notion of rotation number of a homeomorphism of the unit circle dates
back to Poincare. Its properties are derived in Devaney's book \cite{dev}
formally for diffeomorphisms but the proofs are valid for homeomorphisms.
See also \cite{franks} for a helpful survey.

\begin{theorem}
\label{rotnum}Suppose $g:S^{1}\rightarrow S^{1}$ is an orientation
preserving homeomorphism, $\Pi :R\rightarrow S^{1}$ is the covering map $\Pi
(\theta )=e^{2\pi i\theta },x\in R$ and $G_{1},G_{2}:R\rightarrow R$ are
homeomorphisms such that $\Pi (G_{i})=g\Pi .$ Then the following limits
exist and differ by an integer: $\lim_{n\rightarrow \infty }\frac{%
G_{1}^{n}(x)}{n}$ and $\lim_{n\rightarrow \infty }\frac{G_{2}^{n}(x)}{n}.$
This number $(\func{mod}$ 1$)$ is the \textbf{rotation number} of $g$ and is
invariant under the choice of $x.$ If $g$ reverses orientation, then $g$ has
two fixed points and we declare $g$ to have rotation number $0.$ The
rotation number of $g$ is invariant under topological conjugacy ( if $%
h:S^{1}\rightarrow S^{1}$ is a homeomorphism and $g^{\symbol{94}}=hgh^{-1}$
then $rot(g^{\symbol{94}})=rot(g)$). Finally, $rot$ is a continuous function
on the space of homeomorphisms of $S^{1}.$
\end{theorem}

\section{Sets with circular order and rotation number}

Given $4$ distinct point $\{x_{1},x_{2},x_{3},x_{4}\}\subset S^{1}$ declare $%
x_{1}<x_{2}<x_{3}<x_{4}$ if there exists a homeomorphism $g:S^{1}\rightarrow
S^{1}$ such that $g(i^{n})=x_{n}.$

Note it is allowed that $g$\ reverses orientation. Suppose $%
\{J_{1},J_{2},J_{3},J_{4}\}$ is a collection of distinct pairwise disjoint
nonempty subsets of $S^{1}.$ Declare $J_{1}<J_{2}<J_{3}<J_{4}$ if $%
x_{1}<x_{2}<x_{3}<x_{4}$ whenever $x_{i}\in J_{i}.$\bigskip\ 

Suppose $J$ is a proper connected subset of $S^{1}.$ Define $int(J)=J$ if $%
J=\{x\}$ and define $int(J)$ to be the largest open interval contained in $J$
otherwise.

\begin{definition}
Suppose $A^{\symbol{94}}$ is a collection pairwise disjoint subsets of $%
S^{1} $ and $h:A^{\symbol{94}}\rightarrow A^{\symbol{94}}$ is a bijection.
Then $h$ is order preserving if there exists a homeomorphism $%
g:S^{1}\rightarrow S^{1} $ such that $g(J)=h(J)$ for each element $J\in A^{%
\symbol{94}}.$ The homeomorphism $g$ is said to be \textbf{compatible} with $%
h.$ If $\left| A^{\symbol{94}}\right| \geq 3$ then declare $h$ orientation
preserving/reversing iff $g$ is orientation preserving/reversing.
\end{definition}

\begin{lemma}
\label{dokee}Suppose $A^{\symbol{94}}$ is a collection pairwise disjoint
subsets of $S^{1}$ such that $\left| A^{\symbol{94}}\right| \geq 4$ and $%
\cup _{J\in A^{\symbol{94}}}J$ is dense in $S^{1}.$ Suppose $J=int(J)$ for
each $J\in A^{\symbol{94}}.$ Suppose $h:A^{\symbol{94}}\rightarrow A^{%
\symbol{94}}$ is a bijection leaving invariant the set of nontrivial
elements of $A^{\symbol{94}}.$ Then $h$ is order preserving if and only if $%
h(J_{3})$ and $h(J_{4})$ belong to the same component of $S^{1}\backslash
\{h(J_{1})\cup h(J_{2})\}$ whenever $J_{1}<J_{2}<J_{3}<J_{4}.$
\end{lemma}

\begin{proof}
Suppose $h$ is order preserving. Let $g:S^{1}\rightarrow S^{1}$ be a
homeomorphism such that $g(J)=h(J)$ for each $J\in A^{\symbol{94}}.$ Suppose 
$J_{1}<J_{2}<J_{3}<J_{4}.$ Then $J_{3}$ and $J_{4}$ belong to the same
component of $S^{1}\backslash \{J_{1}\cup J_{2}\}.$ Since $g$ is a
homeomorphism $g(J_{3})$ and $g(J_{4})$ belong to the same component of $%
S^{1}\backslash \{g(J_{1})\cup g(J_{2})\}.$ The conclusion follows since $%
g(J_{i})=h(J_{i}).$ Conversely suppose $h(J_{3})$ and $h(J_{4})$ belong to
the same component of $S^{1}\backslash \{h(J_{1})\cup h(J_{2})\}$ whenever $%
J_{1}<J_{2}<J_{3}<J_{4}.$ Suppose $J_{1}<J_{2}<J_{3}<J_{4}.$ Then either $%
h(J_{1})<h(J_{2})<h(J_{3})<h(J_{4})$ or $%
h(J_{1})<h(J_{2})<h(J_{4})<h(J_{3}). $ Suppose in order to obtain a
contradiction that $h(J_{1})<h(J_{2})<h(J_{4})<h(J_{3}).$ Then $h(J_{2})$
and $h(J_{3})$ lie in opposite components of $S^{1}\backslash \{h(J_{1})\cup
h(J_{4}))\}.$ On the other hand by hypothesis $J_{4}<J_{1}<J_{2}<J_{3}$ and
hence $h(J_{2})$ and $h(J_{3})$ lie in the same component of $%
S^{1}\backslash \{h(J_{1})\cup h(J_{4}))\}$ and we have a contradiction.
Thus $h(J_{1})<h(J_{2})<h(J_{3})<h(J_{4})$ and globally the bijection $h$
either preserves or reverses orientation. Let $g:S^{1}\rightarrow S^{1}$ be
the unique homeomorphism mapping $J$ linearly onto $h(J)$ whenever $J$ is a
nontrivial component of $A^{\symbol{94}}.$
\end{proof}

\begin{lemma}
\label{rotg}Suppose $A^{\symbol{94}}$ is a collection of pairwise disjoint
connected subsets of $S^{1}$ such that $\left| A^{\symbol{94}}\right| \geq 3$
and $\cup _{J\in A^{\symbol{94}}}J$ is dense in $S^{1}.$ Suppose $J=int(J)$
whenever $J\in A^{\symbol{94}}$ and suppose $h:A^{\symbol{94}}\rightarrow A^{%
\symbol{94}}$ is an order preserving bijection as demonstrated by the
compatible orientation preserving homeomorphisms $g_{1}:S^{1}\rightarrow
S^{1}$ and $g_{2}:S^{1}\rightarrow S^{1}.$Then $g_{1}$ and $g_{2}$ have the
same rotation number.
\end{lemma}

\begin{proof}
Let $B\subset S^{1}$ denote the complement of the union of the open
intervals of $A^{\symbol{94}}.$ Then $B\neq \emptyset ,$ $g_{1}(B)=B$ and $%
g_{2|B}=g_{1|B}.$ Thus there exists $b^{\symbol{126}}\in R$ and lifts $%
G_{1}:R\rightarrow R$ and $G_{2}:R\rightarrow R$ respectively of $g_{1}$ and 
$g_{2}$ such that $\pi (b^{\symbol{126}})\in B$ and $G_{1}^{n}(b^{\symbol{126%
}})=G_{2}^{n}(b^{\symbol{126}})$ $\forall n.$ Hence by Theorem \ref{rotnum} $%
g_{1}$ and $g_{2}$ have the same rotation number.
\end{proof}

\section{The order preserving bijection $h:A^{\symbol{94}}\rightarrow A^{%
\symbol{94}}$}

\bigskip Define a surjection $p:D^{2}\times \lbrack -1,1]\rightarrow D^{3}$
via $p(x,y,t)=((1-\left| t\right| )x,(1-\left| t\right| )y,t).$ Note $p$
collapses $D^{2}\times \{-1\}$ and $D^{2}\times \{1\}$ to points and is
otherwise injective.

Note the map $\Phi _{X^{\ast }\times (-1,1)}:X^{\ast }\times
(-1,1)\rightarrow D^{2}\times (-1,1)$ defined such that $\Phi (x,t)=(\phi
(x),t)$ can be continuously extended to an injective map $\Phi :\overline{%
X^{\ast }\times (-1,1)}\hookrightarrow \overline{D^{2}\times (-1,1)}$ such
that $\Phi (\infty )=\infty $ and $\Phi (-\infty )=-\infty .$

Define a metric on $X^{\ast }\times R$ via $d((x,t),(y,t))=\max (d^{\ast
}(x,y),\left| s-t\right| ).$

Define a set $A^{\symbol{94}}\subset 2^{\partial D^{2}}$such that $\beta \in
A^{\symbol{94}}$ iff $\beta =int(\phi (X))$ for some path component $%
X\subset \partial X^{\ast }.$

\begin{theorem}
\label{hstar}The homeomorphism $H:\overline{U}\times R\rightarrow \overline{U%
}\times R$ induces a canonical homeomorphism $H^{\ast }:X^{\ast }\times
R\rightarrow X^{\ast }\times R$ such that $H^{\ast }(\partial X^{\ast
}\times R)=\partial X^{\ast }\times R$ and such that $H^{\ast }$ is
extendable to a homeomorphism of $(X^{\ast }\times R)\cup \{\infty ,-\infty
\}.$
\end{theorem}

\begin{proof}
We first obtain, as follows, an induced homeomorphism $H^{\ast }:X^{\ast
}\times R\rightarrow X^{\ast }\times R$ such that $H^{\ast }(U^{\ast }\times
R)=U^{\ast }\times R$ \ and for each $S$ and $T$ there exists $S^{\symbol{94}%
}$ and $T^{\symbol{94}}$ such $H^{\ast }(X^{\ast }\times \lbrack S,T])\cup
H^{\ast -1}(X^{\ast }\times \lbrack S,T])\subset X^{\ast }\times \lbrack S^{%
\symbol{94}},T^{\symbol{94}}].$ Suppose $\{(z_{n},t_{n})\}$ is Cauchy in $%
U^{\ast }\times R$. Since $\{t_{n}\}$ is Cauchy choose $S$ and $T$ such that 
$\forall n(z_{n},t_{n})\in X^{\ast }\times \lbrack S,T]$. Moreover for each $%
n,m$ we may choose a path $\gamma _{nm}\subset U^{\ast }\times \lbrack S,T]$
such that $\gamma _{nm}$ connects $(z_{n},t_{n})$ to $(z_{m},t_{m})$ and diam%
$(\gamma _{nm})<2d((z_{n},t_{n}),(z_{m},t_{m})).$ Hence $\lim_{n,m%
\rightarrow \infty }$ diam$(\gamma _{nm})=0.$ However, since $\overline{U}%
\times \lbrack S,T]$ is compact, both $H_{\overline{U}\times \lbrack S,T]}$
and $H_{\overline{U}\times \lbrack S,T]}^{-1}$ are uniformly continuous.
Thus $\lim_{n,m\rightarrow \infty }$ diam$(H(\gamma
_{nm}))=0=\lim_{n,m\rightarrow \infty }diam(H^{-1}(\gamma _{nm})).$ Hence $H$
and $H^{-1}$ preserve Cauchy sequences in $U^{\ast }\times R$ and thus are
extendable to maps $H^{\ast }:X^{\ast }\times R\rightarrow X^{\ast }\times R$
and $(H^{-1})^{\ast }:X^{\ast }\times R\rightarrow X^{\ast }\times R$ such
that $H^{\ast }((H^{-1})^{\ast })=(H^{-1})^{\ast }H=id.$ Since $id$ is a
bijection it follows that $H^{\ast }$ and $(H^{-1})^{\ast }$ are bijections
and hence homeomorphisms. By definition $H(U\times R)=H^{-1}(U\times
R)=U\times R.$ Since compact subsets of $\overline{U}\times R$ are bounded
we may choose $S^{\symbol{94}}$ and $T^{\symbol{94}}$ such that the
compactum $H(\overline{U}\times \lbrack S,T])\cup H^{-1}(\overline{U}\times
\lbrack S,T])\subset \overline{U}\times \lbrack S^{\symbol{94}},T^{\symbol{94%
}}].$ Note $\overline{U}\times \lbrack S,T]$ (and hence $H(\overline{U}%
\times \lbrack S,T]))$ separate $\overline{U}\times R$ into two components.
Since $\overline{U}\times (T,\infty )$ is connected $H(\overline{U}\times
(T,\infty ))$ intersects and contains exactly one component of $(\overline{U}%
\times R)\backslash (\overline{U}\times \lbrack S^{\symbol{94}},T^{\symbol{94%
}}])$. Thus if $(z_{n},t_{n})\rightarrow \infty $ then $H(z_{n},t_{n})$
converges either to $\infty $ to $-\infty .$ Thus both $H^{\ast }$ and $%
(H^{\ast })^{-1}$ are extendable to maps of $(X^{\ast }\times R)\cup
\{\infty ,-\infty \}$ and both must be homeomorphisms since $H^{\ast
}(H^{\ast -1})$ is fixes pointwise a dense set.
\end{proof}

\begin{corollary}
\label{hbi}The induced homeomorphism $H^{\ast }:X^{\ast }\times R\rightarrow
X^{\ast }\times R$ induces a canonical homeomorphism $H^{\ast \ast }:%
\overline{X^{\ast }\times (-1,1)}\rightarrow \overline{X^{\ast }\times (-1,1)%
}$ and a canonical bijection $h:A^{\symbol{94}}\rightarrow A^{\symbol{94}}.$
\end{corollary}

\begin{proof}
Let $H^{\ast }:X^{\ast }\times R\rightarrow X^{\ast }\times R$ be induced
from $H:\overline{U}\times R\rightarrow \overline{U}\times R$ as in Theorem 
\ref{hstar}. Let $T:R\rightarrow (-1,1)$ be any homeomorphism. Define a
homeomorphism $H^{\ast \ast }:\overline{X^{\ast }\times (-1,1)}\rightarrow 
\overline{X^{\ast }\times (-1,1)}$ such that $H^{\ast \ast
}(x,T(t))=(y,T(s)) $ if $H^{\ast }(x,t)=(y,s).$ Since $H^{\ast \ast
}(\partial X^{\ast }\times (-1,1))=\partial X^{\ast }\times (-1,1)$ the
homeomorphism $H^{\ast \ast }$ permutes the path components of $\partial
X^{\ast }\times (-1,1).$ Each path component of $\partial X^{\ast }\times
(-1,1)$ is of the form $X\times (-1,1) $ where $X$ is a path component of $%
\partial X^{\ast }.$ By definition $\beta \in A^{\symbol{94}}$ iff there
exists a path component $X\subset \partial X^{\ast }$ such that $\beta
=int\phi (X).$ Thus we define a bijection $h:A^{\symbol{94}}\rightarrow A^{%
\symbol{94}}$ satisfying $h(\beta )=\gamma $ if $\beta =int\phi (X)$ and $%
H^{\ast \ast }(X\times (-1,1))=Y\times (-1,1)$ and $\gamma =int\phi (Y).$
\end{proof}

\begin{lemma}
\bigskip \label{okee} Suppose $D\subset D^{3}$ is a topological disk such
that $int(D)\subset int(D^{3}),$ and $\partial D\subset \partial D^{3}.$
Suppose $\{(0,0,1),(0,0,-1)\}\subset \partial D$ and suppose $\alpha
:[0,1]\rightarrow D^{3}\backslash D$ satisfies $\alpha (t)\in int(D^{3})$
iff $0<t<1.$\ Suppose $\beta _{1}$ and $\beta _{2}$ are disjoint nonempty
connected subsets of $S^{1}.$ Suppose $\partial D\subset p((\beta _{1}\cup
\beta _{2})\times \lbrack -1,1]).$ Suppose $\{\alpha (0),\alpha (1)\}\cap
p((\beta _{1}\cup \beta _{2})\times \lbrack -1,1])=\emptyset .$ Then $\pi
_{1}(p^{-1}(\alpha (0),\alpha (1)))$ is contained in a single component of $%
S^{1}\backslash \{\beta _{1}\cup \beta _{2}\}.$
\end{lemma}

\begin{proof}
The topological disk $D$ separates $D^{3}$ into two components and $\partial
D$ separates $\partial D^{3}$ into two components. Moreover the two
components of $\partial D^{3}\backslash \partial D$ are contained in
distinct components of $D^{3}\backslash D$.Hence, since $\iota m(\alpha )$
is connected and $im(\alpha )\cap D=\emptyset ,$ $\partial \alpha $ belongs
to a single component of $\partial D^{3}\backslash \partial D.$ Moreover $%
\partial D^{3}\backslash p((\beta _{1}\cup \beta _{2})\times \lbrack -1,1])$
contains at most two components. No two of these components are contained in
the same component of $\partial D^{3}\backslash \partial D.$ Thus $\{\alpha
(0),\alpha (1)\}$ is contained in a single component of $\partial
D^{3}\backslash p((\beta _{1}\cup \beta _{2})\times \lbrack -1,1]).$ Hence $%
\pi _{1}(p^{-1}(\alpha (0),\alpha (1)))$ is contained in a single component
of $S^{1}\backslash \{\beta _{1}\cup \beta _{2}\}.$
\end{proof}

\begin{corollary}
\label{hop}If $\left| A^{\symbol{94}}\right| \geq 4$ the bijection $h:A^{%
\symbol{94}}\rightarrow A^{\symbol{94}}$ is order preserving.
\end{corollary}

\begin{proof}
Suppose $\gamma _{1}<\gamma _{2}<\gamma _{3}<\gamma _{4}$ and $\gamma
_{i}\in A^{\symbol{94}}$. Let $\iota nt\phi X_{i}=\gamma _{i}$. Choose 4
points $x_{i}\in X_{i}.$ Let $\lambda \subset X^{\ast }$ be an arc
connecting $x_{1}$ to $x_{2}$ such that $int(\lambda )\subset U.$ Let $%
\mathcal{A}\subset X^{\ast }$ be an arc connecting $x_{3}$ to $x_{4}$ such
that $int\mathcal{(A)}\subset U$ and such that $\mathcal{A}\cap \lambda
=\emptyset .$ Let $\Lambda =\{\infty ,-\infty \}\cup (\lambda \times
(-1,1))\subset \overline{X^{\ast }\times (-1,1)}.$ Observe $\Lambda $ is
compact and homeomorphic to $D^{2}$. Let $D=p\Phi (H^{\ast \ast }\Lambda ).$
Let $\alpha =p\Phi (H^{\ast \ast }(\mathcal{A}\times \{0\})).$ Let $\beta
_{i}=h(\gamma _{i}).$ Apply Lemma \ref{okee} to conclude $\beta _{3}$ and $%
\beta _{4}$ belong to the same component of $S^{1}\backslash \{\beta
_{1}\cup \beta _{2}\}.$ Thus by Lemma \ref{dokee} $h$ is order preserving.
\end{proof}

\section{Main results}

Define $G(U)=\{H:\overline{U}\rightarrow \overline{U}|$ $H$ is a
homeomorphism and $H(U)=U\}.$

Define $G(U\times R)=\{H:\overline{U}\times R\rightarrow \overline{U}\times
R|$ $H$ is a homeomorphism and $H(U\times R)=U\times R\}.$

Recalling Theorems \ref{prmend} and \ref{rotnum} define $rot(h):G(U)%
\rightarrow R$ such that $rot(h)$ is the rotation number of $\overline{\psi
h\psi ^{-1}}_{S^{1}}.$

\subsection{The case $\partial U$ is not locally connected}

Suppose $\partial U$ is not locally connected.

Recall the bijection $h:A^{\symbol{94}}\rightarrow A^{\symbol{94}}$ from
Corollary \ref{hbi}. Define a function $Rot:G(U\times R)\rightarrow R$ as
follows.

If $\left| A^{\symbol{94}}\right| \geq 3$ recall Corollary \ref{hop}, note $%
h $ is order preserving, and define $Rot(H)=Rot(g)$ where $%
g:S^{1}\rightarrow S^{1}$ is any homeomorphism compatible with $h.$

If $\left| A^{\symbol{94}}\right| \leq 2$ define $Rot(H)=0$ if $H$ is
orientation reversing and fixes $\{\infty ,-\infty \}$ pointwise or if $H$
is orientation preserving and swaps $\{\infty ,-\infty \}.$ Otherwise define 
$Rot(H)$ to be $0$ or $\frac{1}{2}$ determined, respectively, by whether $%
h=ID_{A^{\symbol{94}}}$ or not.

\subsection{The case $\partial U$ is locally connected}

\bigskip The case $\partial U$ is locally connected deserves special
treatment since $\overline{U}\times R$ admits a homeomorphism extension
theorem analogous to Theorem \ref{prmend} and seems to provide a useful
model of \ `3 dimensional prime end theory' in the sense of \cite{brech}.

\begin{theorem}
\label{pelc}Suppose $\partial U$ is locally connected and $H\in G(U\times
R). $ Then $\Psi H\Psi ^{-1}:int(D^{2})\times R\rightarrow int(D^{2})\times
R $ can be extended to a canonical homeomorphism of the 3 cell $(D^{2}\times
R)\cup \{\infty ,-\infty \}.$
\end{theorem}

\begin{proof}
Since $\partial U$ is locally connected the conformal map $\psi
^{-1}:int(D^{2})\rightarrow U$ can be continuously extended to a surjective
map $\overline{\psi ^{-1}}:D^{2}\rightarrow \overline{U}$ (Theorem 2.1 p20 
\cite{pom}). To show $X^{\ast }$ is compact, it suffices, since $X^{\ast }$
is complete, to show that each sequence in $X^{\ast }$ has a Cauchy
subsequence. Suppose $y_{n}$ is a sequence in $X^{\ast }.$ Let $x_{n}=\phi
(y_{n})$. Let $\{x_{n_{k}}\}\subset D^{2}$ be a Cauchy subsequence of $%
\{x_{n}\}.$ Consider the chords $[x_{n_{l}},x_{n_{k}}]\subset D^{2}.$ Since $%
\overline{\psi ^{-1}}$ is uniformly continuous, diam $(\overline{\psi ^{-1}}[%
x_{n_{l}},x_{n_{k}}])\rightarrow 0.$\ Hence $d^{\ast
}(y_{n_{l}},y_{n_{k}})\rightarrow 0.$ Consequently $X^{\ast }$ is compact
and the injective map from Lemma \ref{smallmain} $\phi :X^{\ast }\rightarrow
D^{2}$ is a homeomorphism. Thus the homeomorphism from Theorem \ref{hstar} $%
H^{\ast }:X^{\ast }\times R\rightarrow X^{\ast }\times R$ induces a
homeomorphism of $D^{2}\times R$ leaving invariant the ends $\{\infty \}$
and $\{-\infty \}.$

To see that $H$ is canonical let $j:U^{\ast }\times R\rightarrow U\times R$
denote identity. Suppose $H_{N}\rightarrow H$ in $G(U\times R).$ Note $%
\{j^{-1}H_{n}j\}$ is uniformly equicontinuous and converges pointwise to $%
j^{-1}Hj$. Hence $H_{n}^{\ast }\rightarrow H^{\ast }$ uniformly. Thus $\Psi
H_{n}^{\ast }\Psi ^{-1}\rightarrow \Psi H^{\ast }\Psi ^{-1}$ uniformly.
\end{proof}

To define rotation number when $\partial U$ is locally connected we seek a
nonempty proper totally disconnected set $B^{\symbol{94}}\subset S^{1}$ such
that $B^{\symbol{94}}\times R$ is invariant under the induced homeomorphism $%
\overline{\Psi H\Psi ^{-1}}:D^{2}\times R\rightarrow D^{2}\times R.$

The point $x\in \partial U$ is a \textbf{cutpoint} if $\partial U\backslash
\{x\}$ not connected. Let $f:D^{2}\rightarrow \overline{U}$ denote the
continuous extension of $\psi ^{-1}.$ The map $f$ is \textbf{light} since
the prime ends of $U$ are in bijective correspondence with points of $S^{1}$
(Theorem 2.15 \cite{pom})$.$ Let $B=\{x\in S^{1}|f(x)$ is a cutpoint of $%
\partial U\}.$

\begin{lemma}
\label{cutpoint}$z\in B$ if and only if $\left| f^{-1}(z)\right| \geq 2.$
\end{lemma}

\begin{proof}
If $z\in \partial U$ and $\left| f^{-1}(z)\right| =1$ then $z\notin B$ is
since $S^{1}\backslash f^{-1}(z)$ is connected and hence $f(S^{1}\backslash
f^{-1}(z))=\partial U\backslash z$ is connected.
\end{proof}

Suppose $\left| f^{-1}(z)\right| \geq 2$ choose $x\neq y$ such that $%
f(y)=f(x)=z.$ Let $[x,y]\subset D^{2}$ denote the chord from $x$ to $y.$ If $%
V$ is a component of $int(D^{2})\backslash \lbrack x,y]$ then $f(\partial
V)\subset \partial (f(V)).$ Thus, since $f$ is light, each complementary
domain of the simple closed curve $f([x,y])$ has nonempty intersection with $%
\partial U.$ Hence $z\in B.$

\begin{lemma}
\label{noon}Suppose $x\neq y,\{x,y\}\subset S^{1}$ and $z=f(x)=f(y).$
Suppose $a$ and $b$ lie in distinct components of $S^{1}\backslash \{x,y\}.$
Suppose $\{a,b\}\cap f^{-1}(z)=\emptyset .$ Then $f(a)\neq f(b).$
\end{lemma}

\begin{proof}
The points $a$ and $b$ belong to the boundaries of distinct components of $%
int(D^{2})\backslash \lbrack x,y].$ These components in turn map into
distinct complementary domains of the simple closed curve $f([x,y]),$ since
points of $int(D^{2})$ can approach the chord $[x,y]$ from distinct sides.
By hypothesis $\{f(a),f(b)\}\cap f([x,y])=\emptyset $ $.$ Thus $f(a)$ and $%
f(b)$ belong to distinct components of $f([x,y])$ and in particular $%
f(a)\neq f(b).$
\end{proof}

\begin{lemma}
If $\partial U$ is not a simple closed curve then $B\neq \emptyset $ and $%
B\neq S^{1}.$
\end{lemma}

\begin{proof}
If $\partial U$ is not a simple closed curve then by Lemma \ref{cutpoint} $f$
is not one to one and hence $B\neq \emptyset .$ Let $\alpha _{0}\subset
S^{1} $ be any closed interval with distinct endpoints $x_{0}$ and $y_{0}$
such that $z_{0}=f(x_{0})=f(y_{0}).$ Proceeding recursively , if possible
let $\alpha _{n+1}\subset \alpha _{n}\backslash f^{-1}(z_{n})$ be a closed
interval with distinct endpoints $x_{n+1}$ and $y_{n+1}$ such that $%
f(x_{n+1})=f(y_{n+1})$ and $\left| x_{n+1}-y_{n+1}\right| <\frac{1}{n+1}.$
Suppose the process never terminates. Let $a=\lim x_{n}=\lim y_{n}.$ Suppose 
$f(b)=f(a).$ By Lemma \ref{noon} $b\in \cap \alpha _{n}.$ Thus $b=a$ and by
Lemma \ref{cutpoint} $a\notin B.$ If the process terminates then there must
exists a nontrivial open interval $(x,y)\subset \alpha _{0}$ such that $%
f_{J} $ is one to one and $f(x)=f(y).$ Let $a\in (x,y)$ and apply Lemma \ref
{noon} to conclude $a\notin B.$
\end{proof}

\begin{lemma}
The following are equivalent.1) $b\notin B.$ 2) Suppose $g:S^{1}\rightarrow 
\overline{U}\times R$ satisfies $g(1)=(f(b),t^{\symbol{94}})$ and $g(\theta
)\in U\times R$ if $\theta \neq 1.$ Then there exists a map $%
G:D^{2}\rightarrow \overline{U}\times R$ such that $G(z)\in U\times R$
whenever $z\in int(D^{2})$ and $G_{S^{1}}=g.$
\end{lemma}

\begin{proof}
Let $\pi _{1}:\overline{U}\times R\rightarrow \overline{U}$ and $\pi _{2}:%
\overline{U}\times R\rightarrow R$ denote the projection maps. Suppose $%
b\notin B$ and $g:S^{1}\rightarrow \overline{U}\times R$ satisfies $\pi
_{1}g(1)=f(b)$ and $g(\theta )\in U\times R$ if $\theta \neq 1.$ By Lemma 
\ref{cutpoint} $\left| f^{-1}(f(b))\right| =1$ and hence the formula $\alpha
=f^{-1}(\pi _{1}g):S^{1}\rightarrow D^{2}$ determines a map. Let $\beta
:D^{2}\rightarrow D^{2}$ be a continuous extension of $\alpha $ such that $%
\beta (z)\in int(D^{2})$ if $z\in int(D^{2}).$ Define $G^{\symbol{94}%
}:D^{2}\rightarrow \overline{U}$ via $G^{\symbol{94}}=f(\beta ).$ Since $R$
is contractible extend $\pi _{2}(g):S^{1}\rightarrow R$ to a map $%
h:D^{2}\rightarrow R.$ Define $G:D^{2}\rightarrow \overline{U}\times R$ via $%
G(z)=(G^{\symbol{94}}(z),h(z)).$

Conversely suppose $b\in B.$ By Lemma \ref{cutpoint} choose $b^{\symbol{94}%
}\in S^{1}$ such that $b^{\symbol{94}}\neq b$ and $f(b^{\symbol{94}})=f(b).$
Let $[b,b^{\symbol{94}}]\subset D^{2}$ denote the chord. Let $g^{\symbol{94}%
} $ map $S^{1}$ homeomorphically onto the simple closed curve $f([b,b^{%
\symbol{94}}])\subset \overline{U}$ such that $g^{\symbol{94}}(1)=f(b)=f(b^{%
\symbol{94}}).$ Both components of $R^{2}\backslash im(g^{\symbol{94}})$
contain points of $\partial U.$ Hence $g^{\symbol{94}}$ is essential in $%
\{f(b)\}\cup (R^{2}\backslash \partial U).$ In particular $G$ cannot exist
since $\pi _{1}(G)$ would show that $g^{\symbol{94}}$ is inessential.
\end{proof}

Since condition 2) is a topological property $B\times R$ maps onto itself
under the induced homeomorphism $\overline{\Psi H\Psi ^{-1}}:D^{2}\times
R\rightarrow D^{2}\times R.$ It follows that if $B^{\symbol{94}}=\{x\in
S^{1}|$ $\{x\}$ is a component of $B$ or $x$ is an endpoint of a nontrivial
interval of $B\}$ then $B^{\symbol{94}}\neq \emptyset $, $B^{\symbol{94}}$
is totally disconnected, and $\overline{\Psi H\Psi ^{-1}}(B^{\symbol{94}%
}\times R)=B^{\symbol{94}}\times R.$ If $\left| B^{\symbol{94}}\right| \neq
2 $ let $g:S^{1}\rightarrow S^{1}$ be any homeomorphism such that $g(x)=y$
whenever $x$ and $y$ are components of $B^{\symbol{94}}$ such that $%
\overline{\Psi H\Psi ^{-1}}(\{x\}\times R)=\{y\}\times R.$ Define $%
Rot(H)=rot(g).$ If $\left| B^{\symbol{94}}\right| =2$ define $Rot(H)=0$ if $%
\overline{\Psi H\Psi ^{-1}}$ is orientation reversing and fixes $\{\infty
,-\infty \}$ pointwise or if $\overline{\Psi H\Psi ^{-1}}$ is orientation
preserving and swaps $\{\infty ,-\infty \}.$ Otherwise define $Rot(H)$ to be 
$0$ or $\frac{1}{2}$ determined, respectively, by whether $\overline{\Psi
H\Psi ^{-1}}$ interchanges the lines of $B^{\symbol{94}}\times R$ or not.

\subsection{Main Theorem on rotation number over $U\times R$}

\begin{theorem}
\label{mainthm}Suppose $\partial U$ is not a simple closed curve. The
function $Rot:G(\overline{U}\times R)\rightarrow R$ is continuous. If $%
F,H\in G(\overline{U}\times R)$ then $Rot(H)=Rot(FHF^{-1})$. If there exists 
$h\in G(\overline{U})$ such that $\forall z,t$ $H(z,t)=(h(z),t)$ then $%
Rot(H)=rot(h).$
\end{theorem}

\begin{proof}
If $\left| A^{\symbol{94}}\right| \geq 3$ then by Corollary \ref{hop} $Rot$
is well defined. Suppose $g,f:S^{1}\rightarrow S^{1}$ are homeomorphisms
compatible, respectively with $H$ and $F.$ Then $fgf^{-1}$ is compatible
with $FGF^{-1}.$ Hence $Rot(H)=rot(g)=rot(fgf^{-1})=Rot(FGF^{-1}).$ For
continuity of $Rot$ suppose $H_{n}\rightarrow H.$ Let $h_{n}:A^{\symbol{94}%
}\rightarrow A^{\symbol{94}}$ and $h:A^{\symbol{94}}\rightarrow A^{\symbol{94%
}}$ denote the corresponding order preserving bijections. If $J$ is a
nontrivial component of $A^{\symbol{94}}$ then $h_{n}(J)$ is eventually
constant. Construct compatible homeomorphisms $g_{n}:S^{1}\rightarrow S^{1}$
and $g:S^{1}\rightarrow S^{1}$ as in Lemma \ref{dokee}. Note eventually $%
g_{n|J}=g_{|J}$ and hence $rot(g_{n})=rot(g)$ eventually. If each component $%
J$ of $A^{\symbol{94}}$ is trivial then $g_{n}\rightarrow g$ pointwise on a
dense set $A^{\symbol{94}}\subset S^{1}$ and hence $g_{n}\rightarrow g$
uniformly. Thus by Theorem \ref{rotnum} $rot(g_{n})\rightarrow rot(g).$

If $\partial U$ is not locally connected and $\left| A^{\symbol{94}}\right|
\leq 2$ the theorem follows from the definition of $Rot.$

If $\partial U$ is locally connected the homeomorphism $\overline{\Psi H\Psi
^{-1}}$ varies continuously with $H.$ The theorem follows from the
definition of $Rot$ and from Theorem \ref{rotnum}.
\end{proof}

\bigskip 

\end{document}